\documentclass[11pt]{amsart}
\usepackage{amssymb}
\usepackage{graphicx}
\usepackage[all]{xy}

\title[Non-deterministic branching and merging]{Homological properties of non-deterministic branchings and mergings in higher dimensional automata}
\author[P. Gaucher]{Philippe Gaucher}
\address{Preuves Programmes et Syst{\`e}mes\\ Universit{\'e} Paris 7--Denis Diderot\\
Case 7014\\2 Place Jussieu\\ 75251 PARIS Cedex 05\\ France}
\email{gaucher@pps.jussieu.fr}
\urladdr{http://www.pps.jussieu.fr/{\~{}}gaucher/}

\newcommand{\C}{\mathcal{C}}
\newcommand{\D}{\mathcal{D}}

\newcommand{\Z}{\mathbb{Z}}

\newcommand{\R}{\mathbb{R}}

\newcommand{\p}\times
\renewcommand{\vec}{\overrightarrow}
\renewcommand{\P}{\mathbb{P}}

\newcommand{\be}{\begin{equation}}
\newcommand{\ee}{\end{equation}}
\newcommand{\bea}{\begin{eqnarray}}
\newcommand{\eea}{\end{eqnarray}}
\newcommand{\beas}{\begin{eqnarray*}}
\newcommand{\eeas}{\end{eqnarray*}}

\newtheorem{thm}{Theorem}[section]
\newtheorem*{thmN}{Theorem}

\newtheorem{prop}[thm]{Proposition}
\newtheorem{lem}[thm]{Lemma}

\newtheorem{cor}[thm]{Corollary}
\newtheorem{rem}[thm]{Remark}

\newtheorem{defn}[thm]{Definition}
\newtheorem{propdef}[thm]{Proposition and Definition}

\newtheorem{nota}[thm]{Notation}
\newcommand{\bd}{\begin{defn}}
\newcommand{\ed}{\end{defn}}
\newcommand{\bcd}{\begin{defn}}
\newcommand{\ecd}{\end{defn}}
\newcommand{\bex}{\begin{exmp}}
\newcommand{\eex}{\end{exmp}}
\newcommand{\bp}{\begin{prop}}
\newcommand{\ep}{\end{prop}}
\newcommand{\bth}{\begin{thm}}
\renewcommand{\eth}{\end{thm}}
\newcommand{\br}{\begin{rem}}
\newcommand{\er}{\end{rem}}
\newcommand{\bpf}{\begin{proof}}
\newcommand{\epf}{\end{proof}}

\newcommand{\fl}[1]{\ar@{->}[l]_{#1}}
\newcommand{\fr}[1]{\ar@{->}[r]^-{#1}}
\newcommand{\fd}[1]{\ar@{->}[d]_{#1}}
\newcommand{\fu}[1]{\ar@{->}[u]^{#1}}
\newcommand{\f}[2]{\ar@{->}[#1]|{#2}}
\newcommand{\ff}[2]{\ar@2{->}[#1]|{#2}}
\newcommand{\frr}[1]{\ar@{->}[rr]^{#1}}

\renewcommand{\top}{{\mathbf{Top}}}

\newcommand{\iso}{\cong}
\newcommand{\lp}{\left(}
\newcommand{\rp}{\right)}

\newcommand{\ot}{\otimes}

\newcommand{\vI}{\vec{I}}

\renewcommand{\geq}{\geqslant}

\newcommand{\Rm}{\mathcal{R}^-}

\newcommand{\cattop}{{\brm{1Cat^{top}_1}}}

\newcommand{\cattopn}{{\brm{1Cat^{top}_1}}}

\newcommand{\homcat}{[\cattopn]}

\def\cartesien{%
  \ar@{-}[]+R+<6pt,-2pt>;[]+RD+<6pt,-6pt>%
  \ar@{-}[]+D+<2pt,-6pt>;[]+RD+<6pt,-6pt>%
}
\def\cocartesien{%
  \ar@{-}[]+L+<-6pt,+2pt>;[]+LU+<-6pt,+6pt>%
  \ar@{-}[]+U+<-2pt,+6pt>;[]+LU+<-6pt,+6pt>%
}

\newcommand{\brm}[1]{\rm{\mathbf{#1}}}

\renewcommand{\top}{{\brm{Top}}}

\newcommand{\dtop}{{\brm{Flow}}}

\newcommand{\tdtop}{{\brm{FLOW}}}
\newcommand{\ttop}{{\brm{TOP}}}

\newcommand{\glob}{{\rm{Glob}}}

\DeclareMathOperator{\sing}{Sing}

\newcommand{\liminj}{\varinjlim}

\makeatletter

\def\varholim@#1#2{%
  \vtop{\m@th\ialign{##\cr
    \hfil$#1\operator@font holim$\hfil\cr
    \noalign{\nointerlineskip\kern1.5\ex@}#2\cr
    \noalign{\nointerlineskip\kern-\ex@}\cr}}%
}
\def\holimproj{%
  \mathop{\mathpalette\varholim@{\leftarrowfill@\textstyle}}\nmlimits@
}
\def\holiminj{%
  \mathop{\mathpalette\varholim@{\rightarrowfill@\textstyle}}\nmlimits@
}

\newskip\@bigflushglue \@bigflushglue = -100pt plus 1fil

\def\bigcentering{\let\\\@centercr\rightskip\@bigflushglue%
\leftskip\@bigflushglue
\parindent\z@\parfillskip\z@skip}

\makeatother

\DeclareMathOperator{\hop}{ho\P} \DeclareMathOperator{\nat}{Nat}
\DeclareMathOperator{\id}{Id}

\DeclareMathOperator{\im}{im}

\begin{document}

\begin{abstract}
The branching (resp. merging) space functor of a flow is a left
Quillen functor. The associated derived functor allows to define the
branching (resp. merging) homology of a flow. It is then proved that
this homology theory is a dihomotopy invariant and that higher
dimensional branchings (resp. mergings) satisfy a long exact sequence.
\end{abstract}

\maketitle

\tableofcontents

\section{Introduction}

The category of \textit{flows} \cite{model3} is an algebraic
topological model of \textit{higher dimensional automata}
\cite{pratt} \cite{rvg}.  Two kinds of mathematical problems are particularly of
importance for such objects: 1) reducing the size of the category of
flows by the introduction of a class of \textit{dihomotopy
equivalences} identifying flows having the same computer-scientific
properties ; 2) investigating the mathematical properties of these
dihomotopy equivalences for instance by constructing related model
category structures and algebraic invariants. For other examples of
similar investigations with different algebraic topological models of
concurrency, cf. for example \cite{mg} \cite{bubenik} \cite{survol}.

This paper is concerned with the second kind of mathematical
problems. Indeed, the purpose of this work is the construction of two
dihomotopy invariants, the \textit{branching homology} $H_*^-(X)$ and
the \textit{merging homology} $H_*^+(X)$ of a flow $X$, detecting the
non-deterministic branching areas (resp. merging areas) of
non-constant execution paths in the higher dimensional automaton
modelled by the flow $X$. Dihomotopy invariance means in the framework
of flows invariant with respect to \textit{weak S-homotopy}
(Corollary~\ref{inv1}) and with respect to \textit{T-homotopy}
(Proposition~\ref{inv2}).

The core of the paper is focused on the case of branchings. The case
of mergings is similar and is postponed to Appendix~\ref{casemerging}.

The \textit{branching space} of a flow is introduced in
Section~\ref{branchingspace} after some reminders about flows
themselves in Section~\ref{rf}. Loosely speaking, the branching space
of a flow is the space of germs of non-constant execution paths
beginning in the same way.  This functor is the main ingredient in the
construction of the branching homology.

However it is badly behaved with respect to \textit{weak S-homotopy
equivalences}, as proved in Section~\ref{bad}.  Therefore it cannot be
directly used for the construction of a dihomotopy invariant. This
problem is overcome in Section~\ref{homotopybranching} by introducing
the \textit{homotopy branching space} of a flow: compare
Theorem~\ref{counter} and Corollary~\ref{inv}. The link between the
homotopy branching space and the branching space is that they coincide
up to homotopy for cofibrant flows, and the latter are the only
interesting and real examples (Proposition~\ref{meme}).

Using this new functor, the branching homology is
finally constructed in Section~\ref{branchinghomology} and it is
proved in the same section and in Section~\ref{verifT} that it is a
dihomotopy invariant (Corollary~\ref{inv1} and
Proposition~\ref{inv2}).

Section~\ref{exbranch} uses the previous construction to establish the
following long exact sequence for higher dimensional branchings:
\begin{thmN}
For any morphism of flows $f:X\longrightarrow Y$, one has the long exact sequence
\beas
&& \dots \rightarrow H_{n}^-(X) \rightarrow H_{n}^-(Y) \rightarrow H_{n}^-(Cf)\rightarrow  \dots \\
&& \dots \rightarrow H_{3}^-(X) \rightarrow H_{3}^-(Y) \rightarrow H_{3}^-(Cf)\rightarrow \\
&& H_{2}^-(X) \rightarrow H_{2}^-(Y) \rightarrow H_{2}^-(Cf)\rightarrow \\
&& H_0(\hop^-X) \rightarrow H_0(\hop^-Y)\rightarrow  H_0(\hop^- Cf)\rightarrow 0.
\eeas
where $Cf$ is the cone of $f$ and where $H_0(\hop^- Z)$ is the free
abelian group generated by the path-connected components of the
homotopy branching space of the flow $Z$.
\end{thmN}
By now, this homological result does not have any known computer
scientific interpretation. But it sheds some light on the potential of
an algebraic topological approach of concurrency.

At last, Section~\ref{example} then gives several examples of
calculation which illustrate the mathematical notions presented here.

Appendix~\ref{correctforShomotopy} is a technical section which proves
that two S-homotopy equivalent flows (which are not necessary
cofibrant) have homotopy equivalent branching spaces. The result is
not useful at all for the core of the paper but is interesting enough
to be presented in an appendix of a paper devoted to branching
homology.

Some familiarity with model categories is required for a good
understanding of this work. However some reminders are included in
this paper. Possible references for model categories are
\cite{MR99h:55031}, \cite{ref_model2} and \cite{MR1361887}. 
The original reference is \cite{MR36:6480}.

\section{The category of flows} 
\label{rf}

In this paper, $\top$ is the category of compactly generated
topological spaces, i.e.  of weak Hausdorff $k$-spaces
(cf. \cite{MR90k:54001},
\cite{MR2000h:55002} and the appendix of \cite{Ref_wH}).

\bd Let $i:A\longrightarrow B$ and $p:X\longrightarrow Y$ be maps in a
category $\C$. Then $i$ has the \textit{left lifting property} (LLP)
with respect to $p$ (or $p$ has the \textit{right lifting property}
(RLP) with respect to $i$) if for any commutative square
\[
\xymatrix{
A\fd{i} \fr{\alpha} & X \fd{p} \\
B \ar@{-->}[ru]^{g}\fr{\beta} & Y}
\]
there exists $g$ making both triangles commutative. \ed

The category $\top$ is equipped with the unique model structure having
the weak homotopy equivalences as weak equivalences and having the
Serre fibrations~\footnote{that is a continuous map having the RLP
with respect to the inclusion $\mathbf{D}^n\p 0\subset
\mathbf{D}^n\p [0,1]$ for any $n\geq 0$ where $\mathbf{D}^n$ is the
$n$-dimensional disk} as fibrations.

\bd\cite{model3}
A {\rm flow} $X$ consists of a compactly generated topological space
$\P X$, a discrete space $X^0$, two continuous maps $s$ and $t$ called
respectively the source map and the target map from $\P X$ to $X^0$
and a continuous and associative map $*:\{(x,y)\in \P X\p \P X;
t(x)=s(y)\}\longrightarrow \P X$ such that $s(x*y)=s(x)$ and
$t(x*y)=t(y)$.  A morphism of flows $f:X\longrightarrow Y$ consists of
a set map $f^0:X^0\longrightarrow Y^0$ together with a continuous map
$\P f:\P X\longrightarrow \P Y$ such that $f(s(x))=s(f(x))$,
$f(t(x))=t(f(x))$ and $f(x*y)=f(x)*f(y)$. The corresponding category
is denoted by $\dtop$. \ed

The topological space $X^0$ is called the \textit{$0$-skeleton} of
$X$. The elements of the $0$-skeleton $X^0$ are called \textit{states}
or \textit{constant execution paths}.  The elements of $\P X$ are
called \textit{non-constant execution paths}. An \textit{initial
state} (resp. a \textit{final state}) is a state which is not the
target (resp. the source) of any non-constant execution path. The
initial flow is denoted by $\varnothing$. The terminal flow is denoted
by $\mathbf{1}$. The initial flow $\varnothing$ is of course the
unique flow such that $\varnothing^0=\P\varnothing=\varnothing$ (the
empty set). The terminal flow is defined by $\mathbf{1}^0=\{0\}$,
$\P\mathbf{1}=\{u\}$ and the composition law $u*u=u$.

\begin{nota} \cite{model3} 
For $\alpha,\beta\in X^0$, let $\P_{\alpha,\beta}X$ be the subspace of
$\P X$ equipped with the Kelleyfication of the relative topology
consisting of the non-constant execution paths $\gamma$ of $X$ with
beginning $s(\gamma)=\alpha$ and with ending $t(\gamma)=\beta$.
\end{nota}

Several examples of flows are presented in Section~\ref{example}. But two 
examples are important for the sequel:

\bd\cite{model3}
Let $Z$ be a topological space. Then the {\rm globe} of $Z$ is the
flow $\glob(Z)$ defined as follows: $\glob(Z)^0=\{0,1\}$,
$\P\glob(Z)=Z$, $s=0$, $t=1$ and the composition law is trivial. The
mapping $\glob:\top\longrightarrow \dtop$ gives rise to a functor in
an obvious way. \ed

\begin{figure}
\begin{center}
\includegraphics[width=7cm]{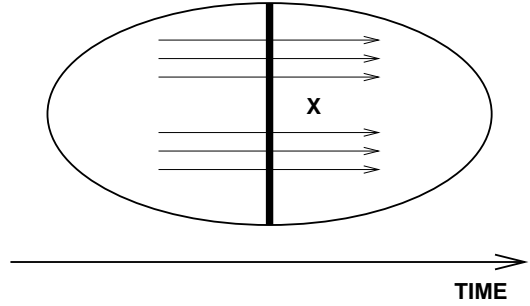}
\end{center}
\caption{Symbolic representation of
$\glob(X)$ for some topological space $X$} 
\label{exglob}
\end{figure}

\begin{nota}\cite{model3} 
If $Z$ and $T$ are two topological spaces, then the flow
\[\glob(Z)*\glob(T)\] is the flow obtained by identifying the final state
of $\glob(Z)$ with the initial state of $\glob(T)$. In other terms,
one has the pushout of flows:
\[
\xymatrix{
\{0\} \fr{0\mapsto 1} \fd{0\mapsto 0} & \glob(Z)\fd{}\\
\glob(T) \fr{} & \cocartesien \glob(Z)*\glob(T)}
\]
\end{nota}

\section{The branching space of a flow} \label{branchingspace}

Loosely speaking, the branching space of a flow is the space of germs
of non-constant execution paths beginning in the same way.

\bp \label{universalpm}
Let $X$ be a flow. There exists a topological space $\P^-X$ unique up
to homeomorphism and a continuous map $h^-:\P X\longrightarrow \P^- X$
satisfying the following universal property:
\begin{enumerate}
\item For any $x$ and $y$ in $\P X$ such that $t(x)=s(y)$, the equality
$h^-(x)=h^-(x*y)$ holds.
\item Let $\phi:\P X\longrightarrow Y$ be a
continuous map such that for any $x$ and $y$ of $\P X$ such that
$t(x)=s(y)$, the equality $\phi(x)=\phi(x*y)$ holds. Then there exists a
unique continuous map $\overline{\phi}:\P^-X\longrightarrow Y$ such that
$\phi=\overline{\phi}\circ h^-$.
\end{enumerate}
Moreover, one has the homeomorphism
\[\P^-X\iso \bigsqcup_{\alpha\in X^0} \P^-_\alpha X\]
where $\P^-_\alpha X:=h^-\lp \bigsqcup_{\beta\in
X^0} \P^-_{\alpha,\beta} X\rp$. The mapping $X\mapsto \P^-X$
yields a functor $\P^-$ from $\dtop$ to $\top$. 
\ep

\bpf
Consider the intersection of all equivalence relations whose graph is
closed in $\P X\p \P X$ and containing the pairs $(x,x*y)$ for any
$x\in\P X$ and any $y\in \P X$ such that $t(x)=s(y)$: one obtains an
equivalence relation $\Rm$. The quotient $\P X/\Rm$ equipped with the
final topology is still a $k$-space since the colimit is the same in
the category of $k$-spaces and in the category of general topological
spaces, and is weak Hausdorff as well since the diagonal of $\P X/\Rm$
is closed in $\P X/\Rm\p \P X/\Rm$.  Let $\phi:\P X\longrightarrow Y$
be a continuous map such that for any $x$ and $y$ of $\P X$ with
$t(x)=s(y)$, the equality $\phi(x)=\phi(x*y)$ holds.  Then the
equivalence relation on $\P X$ defined by ``$x$ equivalent to $y$ if
and only if $\phi(x)=\phi(y)$'' has a closed graph which contains the
graph of $\Rm$. Hence the remaining part of the statement.
\epf

\bd 
Let $X$ be a flow. The topological space $\P^-X$ is called the {\rm
branching space} of the flow $X$. The functor $\P^-$ is called the 
{\rm branching space functor}. \ed

\section{Bad behaviour of the branching space functor}\label{bad}

The purpose of this section is the proof of the following fact:

\bth  \label{counter} There exists a weak S-homotopy equivalence
of flows $f:X\longrightarrow Y$ such that the topological spaces
$\P^-X$ and $\P^-Y$ are not weakly homotopy equivalent. \eth

In other terms, the branching space functor alone is not appropriate
for the construction of dihomotopy invariants.

\begin{lem} Let $Z$ be a flow such that $Z^0=\{\alpha,\beta,\gamma\}$ and
such that $\P Z=\P_{\alpha,\beta}Z\sqcup \P_{\beta,\gamma}Z\sqcup
\P_{\alpha,\gamma}Z$.  Such a flow $Z$ is entirely characterized by
the three topological spaces $\P_{\alpha,\beta}Z$,
$\P_{\beta,\gamma}Z$ and $\P_{\alpha,\gamma}Z$ and the continuous map
$\P_{\alpha,\beta}Z\p\P_{\beta,\gamma}Z\longrightarrow
\P_{\alpha,\gamma}Z$. Moreover, one has the pushout of topological
spaces
\[
\xymatrix{
\P_{\alpha,\beta}Z\p \P_{\beta,\gamma}Z\fr{*} \fd{} & \P_{\alpha,\gamma}Z \fd{}\\
\P_{\alpha,\beta}Z\fr{} & \cocartesien \P_\alpha^-Z}
\]
and the isomorphisms of topological spaces $\P_\beta^-Z\iso
\P_{\beta,\gamma}Z$ and $\P^-Z\iso \P_\alpha^-Z\sqcup \P_\beta^-Z$.
\end{lem}

\bpf It suffices to check that the universal 
property of Proposition~\ref{universalpm} is satisfied by $\P^-Z$.
\epf

For $n\geq 1$, let $\mathbf{D}^n$ be the closed $n$-dimensional disk
and let $\mathbf{S}^{n-1}$ be its boundary. Let
$\mathbf{D}^{0}=\{0\}$. Let $\mathbf{S}^{-1}=\varnothing$ be the empty
space.

\begin{figure}
\begin{center}
\includegraphics[width=10cm]{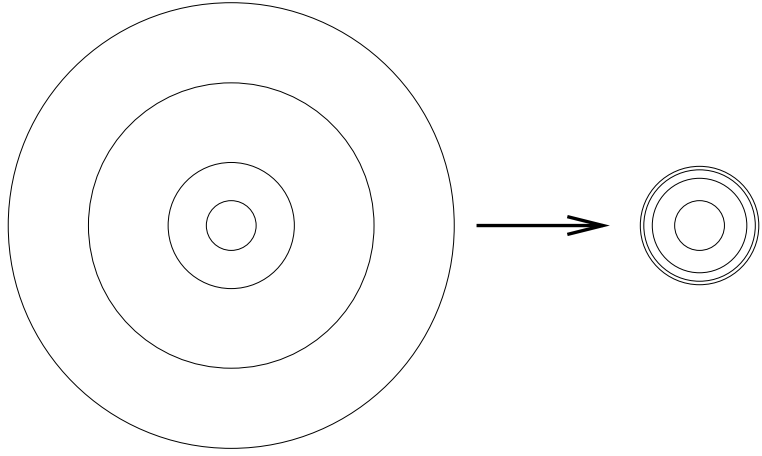}
\end{center}
\caption{$||\phi(x,y)||=\frac{\sqrt{x^2+y^2}}{1+\sqrt{x^2+y^2}}$} 
\label{phi2}
\end{figure}

Let $X$ and $Y$ be the flows defined as follows:
\begin{enumerate}
\item $X^0=Y^0=\{\alpha,\beta,\gamma\}$
\item $\P_{\alpha,\beta}X=\P_{\beta,\gamma}X=\{0\}$
\item $\P_{\alpha,\beta}Y=\P_{\beta,\gamma}Y=\R$
\item $\P_{\alpha,\gamma}X=\P_{\alpha,\gamma}Y=\mathbf{S}^2$
\item the composition law $\P_{\alpha,\beta}X\p \P_{\beta,\gamma}X\longrightarrow \P_{\alpha,\gamma}X$
is given by the constant map $(0,0)\mapsto (0,0,1)\in \mathbf{S}^2$
\item  the composition law $\P_{\alpha,\beta}Y\p \P_{\beta,\gamma}Y\longrightarrow \P_{\alpha,\gamma}Y$ is given by the composite
\[
\xymatrix{
\R\p\R \fr{\phi}&\mathbf{D}^2\backslash \mathbf{S}^1 \ar@{^{(}->}[r]&
\mathbf{D}^2\sqcup_{\mathbf{S}^1} \{(1,0,0)\}\iso \mathbf{S}^2 }
\] where $\phi$ is the homeomorphism (cf. Figure~\ref{phi2}) defined by
\[\phi(x,y)=\lp \frac{x}{1+\sqrt{x^2+y^2}},\frac{y}{1+\sqrt{x^2+y^2}}\rp\]
\end{enumerate}
Then one has the pushouts of compactly generated topological spaces
\[
\xymatrix{\{0\}\p\{0\} \fr{} \fd{} & \mathbf{S}^2\fd{}\\ \{0\}\fr{}& \cocartesien \P_\alpha^-X}
\]
and
\[
\xymatrix{\R\p\R \fr{} \fd{} & \mathbf{S}^2\fd{}\\ \R\fr{}& \cocartesien \P_\alpha^-Y}
\]

\begin{lem}\label{lepushout}
One has the pushout of compactly generated topological spaces
\[
\xymatrix{\R\p\R \fr{} \fd{} & \mathbf{D}^2\sqcup_{\mathbf{S}^1} \{(1,0,0)\}\iso \mathbf{S}^2\fd{}\\ \R\fr{} & \{(1,0,0)\}\cocartesien}
\]
\end{lem}

\bpf 
Let $k\top$ be the category of $k$-spaces. It is well known that the
inclusion functor $i:\top\longrightarrow k\top$ has a left adjoint
$w:k\top\longrightarrow \top$ such that $w\circ i=\id_\top$. So, first
of all, one has to calculate the pushout in the category of
$k$-spaces:
\[
\xymatrix{\R\p\R \fr{} \fd{} & \mathbf{D}^2\sqcup_{\mathbf{S}^1} \{(1,0,0)\}\iso \mathbf{S}^2\fd{}\\ \R\fr{} & X\cocartesien}
\]
and then, one has to prove that $\{(1,0,0)\}\iso w(X)$.

Colimits in $k\top$ are calculated by taking the colimit of the
underlying diagram of sets and by endowing the result with the final
topology. The colimit of the underlying diagram of sets is exactly the
disjoint sum $\R\sqcup\{(1,0,0)\}$. A subset $\Omega$ of
$\R\sqcup\{(1,0,0)\}$ is open for the final topology if and only its
inverse images in $\R$ and $\mathbf{S}^2$ are both open. The inverse
image of $\Omega$ in $\R$ is exactly
$\Omega\backslash\{(1,0,0)\}$. The inverse image of $\Omega$ in $\R\p
\R$ is exactly $\Omega\backslash\{(1,0,0)\}\p \R$. Therefore the
inverse image of $\Omega$ in $\mathbf{S}^2$ is equal to
$\phi(\Omega\backslash\{(1,0,0)\}\p \R)$ if $(1,0,0)\notin \Omega$,
and is equal to $\phi(\Omega\backslash\{(1,0,0)\}\p \R)\cup
\{(1,0,0)\}$ if $(1,0,0)\in \Omega$. 
Therefore, there are now two mutually exclusive cases:
\begin{enumerate}
\item $(1,0,0)\notin \Omega$; in this case, $\Omega$ is open 
if and only if it is open in $\R$
\item $(1,0,0)\in \Omega$; in that case, $\Omega$ is open 
if and only if $\Omega\backslash \{(1,0,0)\}$ is open in $\R$ and 
$\phi(\Omega\backslash\{(1,0,0)\}\p \R)\cup
\{(1,0,0)\}$ is an open of $\mathbf{S}^2$ containing $(1,0,0)$; the 
latter fact is possible if and only if
$\Omega\backslash\{(1,0,0)\}=\R$ (otherwise, if there exists $x\in
\R\backslash (\Omega\backslash\{(1,0,0)\})$, then the straight line
$\phi(\{x\}\p\R)$ tends to $(1,0,0)$ and is not in the inverse image
of $\Omega$).
\end{enumerate}
As conclusion, $X$ is the topological space having the disjoint sum
$\R\sqcup\{(1,0,0)\}$ as underlying set, and a subset $\Omega$ of $X$
is open if and only if $\Omega$ is an open of $\R$ or $\Omega=X$. In
particular, the topological space $X$ is not weak Hausdorff.

Now the topological space $w(X)$ must be determined. It is known that
there exists a natural bijection of sets $\top(w(X),Y)\iso k\top(X,Y)$
for any compactly generated topological space $Y$. Let
$f:X\longrightarrow Y$ be a continuous map. If $Y=\{f((1,0,0))\}$,
then $f$ is a constant map. Otherwise, there exists $y\neq f((1,0,0))$
in $Y$. The singleton $\{y\}$ is closed in $Y$ since the topological
space $Y$ is compactly generated. So $Y\backslash
\{y\}$ is an open of $Y$ containing $f((1,0,0))$. Therefore
$f^{-1}(Y\backslash \{y\})$ is an open of $X$ containing $(1,0,0)$. So
one deduces the equality $f^{-1}(Y\backslash
\{y\})=X$, or equivalently one deduces that $y\notin f(X)$ for any $y\neq
f((1,0,0))$.  This implies again that $f$ is the constant map
$f=f((1,0,0))$. Thus $k\top(X,Y)\iso \top(\{(1,0,0)\},Y)$. The proof
is complete thanks to Yoneda's Lemma.
\epf

\begin{cor} $\P^-X=\mathbf{S}^2\sqcup \{0\}$ and $\P^-Y=\{(1,0,0)\}\sqcup \R$. \end{cor}

\bpf[Proof of Theorem~\ref{counter}] 
It suffices to prove that there exists a weak S-homotopy equivalence
$f$ of flows $X\longrightarrow Y$. Take the identity of
$\{\alpha,\beta,\gamma\}$ on the $0$-skeleton.  Take the identity of
$\mathbf{S}^2$ for the restriction
$f:\P_{\alpha,\gamma}X\longrightarrow\P_{\alpha,\gamma}Y$. Let
$(u,v)\in \R\p\R$ such that $\phi(u,v)=(0,0,1)$. Then it suffices to
put $f(0)=u$ for $0\in\P_{\alpha,\beta}X$ and $f(0)=v$ for $0\in
\P_{\beta,\gamma}X$.
\epf

The reader must not be surprised by the result of this
section. Indeed, the branching space is given by a colimit. And it is
well-known that colimits are badly behaved with respect to weak
equivalences and that they must be replaced by homotopy colimits in
algebraic topology.

\section{The homotopy branching space} \label{homotopybranching}

Let us denote by $Q$ the cofibrant replacement functor of any model
structure.

\bd \cite{MR99h:55031} \cite{ref_model2} \cite{MR1361887}
An object $X$ of a model category $\C$ is {\rm cofibrant}
(resp. {\rm fibrant}) if and only if the canonical morphism
$\varnothing\longrightarrow X$ from the initial object of $\C$ to $X$
(resp. the canonical morphism $X\longrightarrow \mathbf{1}$ from $X$
to the final object $\mathbf{1}$) is a cofibration (resp.  a
fibration). \ed

In particular, in any model category, the canonical morphism
$\varnothing\longrightarrow X$ where $\varnothing$ is the initial
object) functorially factors as a composite
$\varnothing\longrightarrow Q(X) \longrightarrow X$ of a cofibration
$\varnothing\longrightarrow Q(X)$ followed by a trivial fibration
$Q(X) \longrightarrow X$.

\begin{propdef} \cite{MR99h:55031}
\cite{ref_model2} \cite{MR1361887} \label{rappel2}
A Quillen adjunction is a pair of adjoint functors
$F:\C\rightleftarrows \D:G$ between the model categories $\C$ and $\D$
such that one of the following equivalent properties holds:
\begin{enumerate}
\item if $f$ is a cofibration (resp. a trivial cofibration), then so
is $F(f)$
\item if $g$ is a fibration (resp. a trivial fibration), then so
is $G(g)$.
\end{enumerate}
One says that $F$ is a {\rm left Quillen functor}.  One says that $G$
is a {\rm right Quillen functor}. Moreover, any left Quillen functor
preserves weak equivalences between cofibrant objects and any right
Quillen functor preserves weak equivalences between fibrant objects.
\end{propdef}

The fundamental tool of this section is the:

\bth \cite{model3} \label{rappel}
There exists one and only one model structure on $\dtop$ such that 
\begin{enumerate}
\item 
the weak equivalences are the so-called {\rm weak S-homotopy
equivalences}, that is the morphisms of flows $f:X\longrightarrow Y$
such that $f^0:X^0\longrightarrow Y^0$ is a bijection and such that
$\P f:\P X\longrightarrow \P Y$ is a weak homotopy equivalence of
topological spaces
\item the fibrations are the morphisms of flows $f:X\longrightarrow Y$ such that 
$\P f:\P X\longrightarrow \P Y$ is a (Serre) fibration of topological
spaces. 
\end{enumerate}
Any flow is fibrant for this model structure.  
\eth

\bd\label{defS} \cite{model3}
The notion of homotopy between cofibrant-fibrant flows is called
{\rm S-homo\-topy}.
\ed

\bth\label{Pleft} 
The branching space functor $\P^-:\dtop\longrightarrow \top$ is a left
Quillen functor. \eth

\bpf 
One has to prove that there exists a functor $C^-:\top\longrightarrow
\dtop$ such that the pair of functors
$\P^-:\dtop\rightleftarrows\top:C^-$ is a Quillen adjunction.

Let us define the functor
$C^-:\top\longrightarrow \dtop$ as follows:
$C^-(Z)^0=\{0\}$, $\P C^-(Z)=Z$ with the composition law
${\rm pr}_1:(z_1,z_2)\mapsto z_1$. Indeed, one has
${\rm pr}_1({\rm pr}_1(z_1,z_2),z_3)={\rm pr}_1(z_1,{\rm pr}_1(z_2,z_3))=z_1$.

A continuous map $f:\P^-X\longrightarrow Z$ gives rise to a continuous
map $f\circ h^-:\P X\longrightarrow Z$ such that
\[f(h^-(x*y))=f(h^-(x))={\rm pr}_1(f(h^-(x)),f(h^-(y)))\]
which provides the set map
\[\top(\P^-X,Z)\longrightarrow\dtop(X,C^-(Z)).\] Conversely, if
$g\in\dtop(X,C^-(Z))$, then $\P g:\P X\longrightarrow \P C^-(Z)=Z$
satisfies \[\P g(x*y)={\rm pr}_1(\P g(x),\P g(y))=\P g(x).\] Therefore
$\P g$ factors uniquely as a composite $\P X\longrightarrow
\P^-X\longrightarrow Z$ by Proposition~\ref{universalpm}. So
one has the natural isomorphism of sets
\[\top(\P^-X,Z)\iso \dtop(X,C^-(Z)).\]

A morphism of flows $f:X\longrightarrow Y$ is a fibration if and only
if $\P f:\P X\longrightarrow \P Y$ is a fibration by
Theorem~\ref{rappel}. Therefore $C^-$ is a right Quillen functor and
$\P^-$ is a left Quillen functor by Proposition~\ref{rappel2}. \epf

\bd
The {\rm homotopy branching space} $\hop^- X$ of a flow $X$ is by
definition the topological space $\P^-Q(X)$.  If $\alpha\in X^0$, let 
$\hop^-_\alpha X=\P^-_\alpha Q(X)$. \ed

\begin{cor} \label{inv}
Let $f:X\longrightarrow Y$ be a weak S-homotopy equivalence of flows.
Then $\hop^- f:\hop^- X\longrightarrow \hop^- Y$ is a homotopy
equivalence between cofibrant topological spaces. \end{cor}

\bpf The morphism of flows $Q(f)$ is a weak S-homotopy equivalence between 
cofibrant flows. Since $\P^-$ is a left Quillen adjoint, the morphism
$\hop^- f:\hop^- X\longrightarrow \hop^- Y$ is then a weak homotopy
equivalence between cofibrant topological spaces, and therefore a homotopy 
equivalence by Whitehead's theorem. \epf

\begin{cor}\label{preho}
Let $X$ be a diagram of flows. Then there exists an isomorphism of
flows $\liminj \P^-(X)\iso \P^-(\liminj X)$ where $\liminj$ is the
colimit functor and there exists a homotopy equivalence between the
cofibrant topological spaces $\holiminj\hop^-(X)$ and $\hop^-
(\holiminj X)$ where $\holiminj$ is the homotopy colimit functor.
\end{cor}

The reader does not need to know what a general homotopy colimit is
because Corollary~\ref{preho} will be used only for homotopy
pushout. And a definition of the latter is recalled in
Section~\ref{exbranch}.  Corollary~\ref{preho} is the homotopic analog
of the well-known fact of category theory saying that a left adjoint
commutes with any colimit.

\section{Construction of the branching homology and weak S-homotopy}
\label{branchinghomology}

In this section, we construct the branching homology of a flow and we
prove that it is invariant with respect to weak S-homotopy
equivalences (cf. Theorem~\ref{rappel}).

\bd\label{hombrdef} 
Let $X$ be a flow. Then the $(n+1)$-th branching homology group
$H_{n+1}^-(X)$ is defined as the $n$-th homology group of the
augmented simplicial set $\mathcal{N}^-_*(X)$ defined as follows:
\begin{enumerate}
\item $\mathcal{N}^-_n(X)=\sing_n(\hop^-X)$ for $n\geq 0$
\item $\mathcal{N}^-_{-1}(X)=X^0$
\item the augmentation map $\epsilon:\sing_0(\hop^-X)\longrightarrow X^0$
is induced by the mapping $\gamma\mapsto s(\gamma)$ from $\hop^-X=\sing_0(\hop^-X)$
to $X^0$
\end{enumerate}
where $\sing(Z)$ denotes the singular simplicial nerve of a given
topological space $Z$ \cite{MR2001d:55012}. In other terms, 
\begin{enumerate}
\item for $n\geq 1$, $H_{n+1}^-(X):=H_n(\hop^-X)$
\item  $H_1^-(X):=\ker(\epsilon)/\im\lp\partial:\mathcal{N}^-_1(X)\rightarrow
\mathcal{N}^-_0(X)\rp$
\item $H_0^-(X):=\Z(X^0)/\im(\epsilon)$.
\end{enumerate}
where $\partial$ is the simplicial differential map, where $\ker(f)$ is the kernel 
of $f$ and where $\im(f)$ is the image of $f$. 
\ed

\bp For any flow $X$, $H_0^-(X)$ is the free abelian group generated by the final 
states of $X$. \ep

\bpf Obvious. \epf

Let us denote by $\widetilde{H}_*(Z)$ the reduced homology of a
topological space $Z$, that is the homology group of the augmented
simplicial nerve $\sing(Z)\longrightarrow \{0\}$ (cf. for instance 
\cite{Rotman} definition p. 102). Then one has:

\bp For any flow $X$, there exists a natural isomorphism of abelian groups 
\[H_{n+1}^-(X)\iso \bigoplus_{\alpha\in X^0} \widetilde{H}_n(\hop^-_\alpha X)\] 
for any $n\geq 0$.
\ep

\bpf For $n\geq 1$, one has 
\[\bigoplus_{\alpha\in X^0} \widetilde{H}_n(\hop^-_\alpha X)\iso \bigoplus_{\alpha\in X^0} {H}_n(\hop^-_\alpha X)\iso H_n\lp \bigoplus_{\alpha\in X^0} \hop^-_\alpha X\rp\]
hence the result for $n\geq 1$ by Definition~\ref{hombrdef} and the
$X^0$-grading of $\hop^- X$. For $n=0$, this is a straightforward
consequence of Definition~\ref{hombrdef} and of the definition of the
homology of an augmented simplicial set.
\epf

\bp
Let $f:X\longrightarrow Y$ be a weak S-homotopy equivalence of flows. Then 
$\mathcal{N}^-(f):\mathcal{N}^-(X)\longrightarrow \mathcal{N}^-(Y)$ is a 
homotopy equivalence of augmented simplicial nerves. \ep

\bpf 
This is a consequence of Corollary~\ref{inv} and of the fact that the
singular nerve functor is a right Quillen functor. \epf

\begin{cor} \label{inv1}
Let $f:X\longrightarrow Y$ be a weak S-homotopy equivalence of
flows. Then $H_n^-(f):H_n^-(X)\longrightarrow H_n^-(Y)$ is an
isomorphism for any $n\geq 0$. \end{cor}

\section{Branching homology and T-homotopy}\label{verifT}

In this section, we prove that the branching homology is invariant
with respect to T-homotopy equivalences (cf.  Definition~\ref{defT}).

\begin{figure}
\begin{center}
\includegraphics[width=5cm]{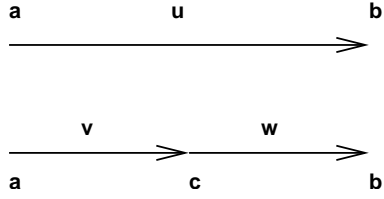}
\end{center}
\caption{Simplest example of T-homotopy equivalence} 
\label{ex1}
\end{figure}

The most elementary example of T-homotopy equivalence which is not
inverted by the model structure of Theorem~\ref{rappel} is the unique
morphism $\phi$ dividing a directed segment in a composition of two
directed segments (Figure~\ref{ex1} and Notation~\ref{phi})

\begin{nota} \label{phi} 
The morphism of flows $\phi:\vI\longrightarrow \vI*\vI$ is the unique
morphism $\phi:\vI\longrightarrow \vI*\vI$ such that
$\phi([0,1])=[0,1]*[0,1]$ where the flow $\vI=\glob(\{[0,1]\})$ is the
directed segment. It corresponds to Figure~\ref{ex1}. \end{nota}

\bd
Let $X$ be a flow. Let $A$ and $B$ be two subsets of $X^0$. One says that
$A$ is {\rm surrounded} by $B$ (in $X$) if for any $\alpha\in A$,
either $\alpha \in B$ or there exists execution paths $\gamma_1$ and
$\gamma_2$ of $\P X$ such that $s(\gamma_1)\in B$,
$t(\gamma_1)=s(\gamma_2)=\alpha$ and $t(\gamma_2)\in B$. We denote
this situation by $A\lll B$.  \ed

\bd\cite{model2}\label{defT}
A morphism of flows $f:X\longrightarrow Y$ is a {\rm T-homotopy
equivalence} if and only if the following conditions are satisfied:
\begin{enumerate}
\item The morphism of flows $f:X\longrightarrow
Y\!\restriction_{(f(X^0)}$ is an isomorphism  of flows. In
particular, the set map $f^0:X^0\longrightarrow Y^0$ is one-to-one.
\item For $\alpha\in Y^0\backslash f(X^0)$, the topological spaces
$\P^-_\alpha Y$ and $\P^+_\alpha Y$ (cf. Proposition~\ref{universalmp}
and Definition~\ref{defplus}) are singletons.
\item $Y^0\lll f(X^0)$.
\end{enumerate}
\ed

\bp\label{inv2} 
Let $f:X\longrightarrow Y$ be a T-homotopy equivalence. Then for any
$n\geq 0$, the linear map $H_n^-(f):H_n^-(X)\longrightarrow H_n^-(Y)$
is an isomorphism.
\ep

\bpf 
For any $\alpha\in X^0$, the continuous map $\hop^-_\alpha
X\longrightarrow \hop^-_\alpha Y$ is a weak homotopy equivalence. So
for $n\geq 1$, one has
\[H_{n+1}^-(X)\iso H_n(\hop^-X)\iso \bigoplus_{\alpha\in X^0} H_n(\hop^-_\alpha X)
\iso \bigoplus_{\alpha\in Y^0} H_n(\hop^-_\alpha Y)\iso H_{n+1}^-(Y)
\]
since for $\alpha\in Y^0\backslash f(X^0)$, the $\Z$-module
$H_n(\hop^-_\alpha Y)$ vanishes.

The augmented simplicial set $\mathcal{N}^-_*(X)$ is clearly
$X^0$-graded. So the branching homology is $X^0$-graded as well. Thus
one has
\[H_1^-(X)=\bigoplus_{\alpha\in X^0} G^\alpha H_1^-(X)\]
with
\[G^\alpha H_1^-(X)\iso \ker\lp \sing_0(\hop^-_\alpha X)\rightarrow \Z\{\alpha\}\rp / \im\lp \Z \sing_1(\hop^-_\alpha X)\rightarrow \Z \sing_0(\hop^-_\alpha X)\rp.\]
So one has the short exact sequences
\[0\rightarrow G^\alpha H_1^-(X)\rightarrow H_0(\hop^-_\alpha X)\rightarrow \Z \hop^-_\alpha X/\ker(s)\rightarrow 0\]
for $\alpha$ running over $X^0$. If $\alpha\in Y^0\backslash f(X^0)$,
then $H_0(\hop^-_\alpha Y)=\Z$.  In this case, $s:\hop^-_\alpha
Y\longrightarrow \{\alpha\}$ so $\Z \hop^-_\alpha Y/\ker(s)\iso \Z$.
Therefore $G^\alpha H_1^-(Y)=0$.

At last, if $\alpha\in Y^0\backslash f(X^0)$, then $\alpha$ belongs to
$\im(s)$ because $Y^0\lll f(X^0)$. Hence the result.
\epf

\begin{cor} The branching homology is a dihomotopy invariant. \end{cor}

\bpf There are two kinds of dihomotopy equivalences in the framework of flows: 
the weak S-homotopy equivalences and the T-homotopy equivalences
\cite{model2}. This corollary is then a consequence of Corollary~\ref{inv1} and 
Proposition~\ref{inv2}. \epf

The reader maybe is wondering why the singular homology of the
homotopy branching space is not taken as definition of the branching
homology.

\bp 
The functor $X\mapsto H_0(\hop^-X)$ is invariant with respect to weak
S-homotopy, but not with respect to T-homotopy equivalences. \ep

\bpf The first part of the statement is a consequence of Corollary~\ref{inv}.  For the
second part of the statement, let us consider the morphism of flows
$\phi:\vI\longrightarrow \vI*\vI$ dividing the directed segment in two
directed segments. Then $H_0(\hop^-\vI)=\Z$ (the path-connected
components of $\P\vI=\{u\}$) and
$H_0(\hop^-(\vI*\vI))=\Z\oplus\Z$ (the path-connected
components of $\P^-(\vI*\vI)=\{v=v*w,w\}$). 
\epf

\section{Long exact sequence for higher dimensional branchings}
\label{exbranch}

\begin{lem} \label{lien}
One has: 
\begin{enumerate}
\item if 
\[
\xymatrix{
U\fr{}\fd{} & V \fd{}\\
W \fr{} & \cocartesien X}
\] 
is a pushout diagram of topological spaces, then 
\[
\xymatrix{
\glob(U)\fr{}\fd{} & \glob(V) \fd{}\\
\glob(W) \fr{} & \cocartesien \glob(X)}
\] 
is a pushout diagram of flows
\item if $g:U\longrightarrow V$ is a cofibration of topological spaces, then 
$\glob(g):\glob(U)\longrightarrow \glob(V)$ is a cofibration of flows 
\item if $U$ is a cofibrant topological space, then $\glob(U)$ is a cofibrant 
flow 
\item there exists a cofibrant replacement functor $Q$ of $\top$ such that 
$Q(\glob(U))=\glob(Q(U))$ for any topological space $U$.  
\end{enumerate}
\end{lem}

\bpf The diagram of sets 
\[
\xymatrix{
\{0,1\}=\glob(U)^0\fr{}\fd{} & \{0,1\}=\glob(V)^0 \fd{}\\
\{0,1\}=\glob(W)^0 \fr{} &  \{0,1\}=\glob(X)^0}
\] 
is a square of constant set maps. Therefore the corresponding pushout
of globes does not create any new non-constant execution paths. Hence
the first assertion.

If $g:U\longrightarrow V$ is a cofibration of topological spaces, then
$g$ is a retract of a transfinite composition of pushouts of morphisms
of $I=\{\mathbf{S}^{n-1}\subset \mathbf{D}^{n},n\geq 0\}$, and
therefore $\glob(g)$ is a retract of a transfinite composition of
pushouts of morphisms of $\{\glob(\mathbf{S}^{n-1})\subset
\glob(\mathbf{D}^{n}),n\geq 0\}$. Since the model structure of Theorem~\ref{rappel} 
is cofibrantly generated with set of
generating cofibrations $I^{gl}_+=\{\glob(\mathbf{S}^{n-1})\subset
\glob(\mathbf{D}^{n}),n\geq 0\}\cup\{R,C\}$ where $R:\{0,1\}\longrightarrow \{0\}$ 
and $C:\varnothing\longrightarrow \{0\}$, the morphism of flows 
$\glob(g):\glob(U)\longrightarrow \glob(V)$ is a cofibration of
flows. Hence the second assertion.

The third assertion is a consequence of the second one and of the fact
that $C:\varnothing\longrightarrow
\{0\}$ is a cofibration.

The cofibrant replacement functor $Q$ of $\dtop$ is obtained by
applying the small object argument for $I^{gl}_+$ with the cardinal
$2^{\aleph_0}$ (\cite{model3} Proposition~11.5). Let $X$ be a
flow. Let $X:2^{\aleph_0}\longrightarrow \dtop$ be the
$2^{\aleph_0}$-sequence with $X^0=\varnothing$ and for any ordinal
$\lambda<2^{\aleph_0}$ by the pushout diagram
\[
\xymatrix{
\bigsqcup_{k\in K} C_k \fr{}\fd{} & X^\lambda\fd{}\\
\bigsqcup_{k\in K} D_k \ar@{=}[d]\fr{} & \cocartesien X^{\lambda+1} \fd{}\\
\bigsqcup_{k\in K} D_k \fr{} & X}
\]
where $K$ is the set of morphisms (i.e. of commutative squares) from a
morphism of $I^{gl}_+$ to the morphism $X^\lambda\longrightarrow
X$. Then $Q(X)=X^{2^{\aleph_0}}$. Pick a topological space $U$ and
consider $X=\glob(U)$. Let $X^0=\varnothing$. Then $X^1=\{0\}\sqcup
\{1\}=\glob(\varnothing)$. Let $U^0=\varnothing$. Let $U:2^{\aleph_0}\longrightarrow
\top$ be the $2^{\aleph_0}$-sequence giving the cofibrant replacement functor of the 
topological space $U$ obtained by applying the small object argument
for $I=\{\mathbf{S}^{n-1}\subset \mathbf{D}^n,n\geq 0\}$ with the
cardinal $2^{\aleph_0}$ (the cardinal $\aleph_0$ is sufficient to
obtain a cofibrant replacement functor in $\top$). Then an easy
transfinite induction proves that $\glob(U^\lambda)=X^{\lambda+1}$. So
$\glob(U^{2^{\aleph_0}})=Q(X)$. The proof of the last assertion is
complete because the functor $U\mapsto U^{2^{\aleph_0}}$ is a
cofibrant replacement functor of $\top$ since $2^{\aleph_0}\geq
\aleph_0$.
\epf

\begin{lem}  (Calculating a homotopy pushout) \label{key} 
In a model category $\mathcal{M}$, the homotopy pushout of the diagram
\[
\xymatrix{
A\fr{i} \fd{} & B \\
C & }
\] 
is homotopy equivalent to the  pushout of the diagram 
\[
\xymatrix{
Q(A)\fr{} \fd{} & Q(B) \\
Q(C) & }
\] 
where $Q$ is a cofibrant replacement functor of $\mathcal{M}$. 
\end{lem}

\bpf Consider the three-object category $\mathcal{B}$
\[
\xymatrix{
1\fr{} \fd{} & 2 \\ 0}
\] 
Let $\mathcal{M}^\mathcal{B}$ be the category of diagrams of objects
of $\mathcal{M}$ based on the category $\mathcal{B}$, or in other
terms the category of functors from $\mathcal{B}$ to
$\mathcal{M}$. There exists a model structure on
$\mathcal{M}^\mathcal{B}$ such that the colimit functor
$\liminj:\mathcal{M}^\mathcal{B}\longrightarrow \mathcal{M}$ is a left
Quillen functor and such that the cofibrant objects are the functors
$F:\mathcal{B}\longrightarrow \C$ such that $F(0)$, $F(1)$ and $F(2)$
are cofibrant in $\C$ and such that $F(1\longrightarrow 2)$ is a
cofibration of $\mathcal{M}$: cf. the proof of the Cube Lemma
\cite{MR99h:55031} \cite{ref_model2}. Hence the result.  \epf

\bd Let $f:X\longrightarrow Y$ be a morphism of flows. The {\rm cone} $Cf$ of $f$
is the homotopy pushout in the category of flows
\[
\xymatrix{
X\fr{f} \fd{} & Y \fd{} \\
\textbf{1} \fr{} & Cf
}
\]
where $\textbf{1}$ is the terminal flow.
\ed

\begin{nota} Let $Z$ be a topological space. Let us denote by
$\underline{L}(Z)$ the pushout
\[
\xymatrix{
\{0,1\} \fd{R} \fr{} & \glob(Z) \fd{} \\
\{0\}  \fr{} & \underline{L}(Z)
}
\]
The $0$-skeleton of $\underline{L}(Z)$ is $\{0\}$ and the path space
of $\underline{L}(Z)$ is $Z\sqcup (Z\p Z)\sqcup (Z\p Z\p Z) \sqcup \dots$.
\end{nota}

\begin{lem} 
Let $g:U\longrightarrow V$ be a cofibration between cofibrant
topological spaces. Then the cone of
$\glob(g):\glob(U)\longrightarrow\glob(V)$ is S-homotopy equivalent to
$\underline{L}(V/U)$.
\end{lem}

\bpf The diagram of flows 
\[
\xymatrix{
Q(\glob(U))\fd{} \fr{} & Q(\glob(V))  \\
Q(\mathbf{1}) & }
\]
induces the diagram of topological spaces 
\[
\xymatrix{
\P Q(\glob(U))\fd{} \fr{} & \P Q(\glob(V))  \\
\P Q(\mathbf{1}) & }
\]
By Lemma~\ref{lien}, one can suppose that $Q(\glob(U))= \glob(Q(U))$
and $Q(\glob(V))= \glob(Q(V))$. Hence one can consider the pushout
diagram of cofibrant topological spaces
\[
\xymatrix{
Q(U) \fd{}\fr{Q(g)} & Q(V) \fd{} \\
\P Q(\mathbf{1}) \fr{} & \cocartesien Z}
\]
By Lemma~\ref{key}, the topological space $Z$ is cofibrant and is
homotopy equivalent to the cone of $g$, that is $V/U$. Since
$Q(\mathbf{1})^0=\{0\}$, one deduces the pushout diagram of flows
\[
\xymatrix{
Q(\glob(U))\fd{} \fr{} & Q(\glob(V)) \fd{} \\
Q(\mathbf{1}) \fr{} & \cocartesien\underline{L}(Z)}
\]
Again by Lemma~\ref{key}, and because $\glob(g)$ is a cofibration of
flows, the flow $\underline{L}(Z)$ is cofibrant and S-homotopy
equivalent to the cone of $\glob(g)$. It then suffices to observe that the flows 
$\underline{L}(Z)$ and $\underline{L}(V/U)$ are S-homotopy equivalent to complete the proof. 
\epf

\begin{lem} \label{point} 
The homotopy branching space of the terminal flow is
contractible. \end{lem}

\bpf Consider the homotopy pushout of flows
\[
\xymatrix{
\glob(U) \fr{\glob(g)} \fd{} & \glob(V) \fd{} \\
\textbf{1} \fr{} & \underline{L}(V/U)
}
\]
where $g:U\longrightarrow V$ is a cofibration between cofibrant
topological spaces. The functor $\hop^-$ preserves homotopy pushouts
by Corollary~\ref{preho}.  Therefore one obtains the homotopy pushout
of topological spaces
\[
\xymatrix{
\hop^-\glob(U) \fr{} \fd{} & \hop^-\glob(V) \fd{} \\
\hop^-\textbf{1} \fr{} & \hop^-\underline{L}(V/U)
}
\]
Since $U$ is cofibrant, $\glob(U)$ is cofibrant as well, therefore
$Q(\glob(U))$ is S-homotopy equivalent to $\glob(U)$. So the space
$\hop^-\glob(U)=\P^- Q(\glob(U))$ is homotopy equivalent to $\P^-
Q(\glob(U))=U$. Since $V/U$ is a cofibrant space as well, the topological 
space 
\[\P \underline{L}(V/U)\iso V/U \sqcup (V/U\p V/U) \sqcup (V/U\p V/U \p V/U) \p \dots \]
is cofibrant as well. So $\hop^-\underline{L}(V/U)$ is homotopy equivalent to
$V/U$. One obtains the homotopy pushout of topological spaces
\[
\xymatrix{
U \fr{g} \fd{} & V \fd{} \\
\hop^-\textbf{1} \fr{} &  V/U
}
\]
for any cofibration $g:U\longrightarrow V$ between cofibrant
spaces. Take for $g$ the identity of $\{0\}$.  One deduces that
$\hop^-\textbf{1}$ is homotopy equivalent to $V/U$, that is to say a
point.
\epf

\begin{lem} \label{commutecone} 
Let $f:X\longrightarrow Y$ be a morphism of flows. Let $Cf$ be the
cone of $f$.  Then the homotopy branching space $\hop^- (Cf)$ of $Cf$
is homotopy equivalent to the cone $C (\hop^-f)$ of
$\hop^-f:\hop^-X\longrightarrow \hop^-Y$.
\end{lem}

\bpf Consider the homotopy pushout of flows
\[
\xymatrix{
X\fr{f} \fd{} & Y \fd{} \\
\textbf{1} \fr{} & Cf
}
\]
Using Corollary~\ref{preho}, one obtains the homotopy pushout of topological
spaces
\[
\xymatrix{
\hop^- X\fr{\hop^-f} \fd{} & \hop^- Y \fd{} \\
\hop^- \textbf{1} \fr{} & \hop^- Cf
}
\]
The proof is complete with Lemma~\ref{point}.
\epf

\bth (Long exact sequence for higher dimensional branchings) 
For any morphism of flows $f:X\longrightarrow Y$, one has the long exact sequence
\beas
&& \dots \rightarrow H_{n}^-(X) \rightarrow H_{n}^-(Y) \rightarrow H_{n}^-(Cf)\rightarrow  \dots \\
&& \dots \rightarrow H_{3}^-(X) \rightarrow H_{3}^-(Y) \rightarrow H_{3}^-(Cf)\rightarrow \\
&& H_{2}^-(X) \rightarrow H_{2}^-(Y) \rightarrow H_{2}^-(Cf)\rightarrow \\
&& H_0(\hop^-X) \rightarrow H_0(\hop^-Y)\rightarrow  H_0(\hop^- Cf)\rightarrow 0.
\eeas
\eth

\bpf If $g:U\rightarrow V$ is a continuous map, then it is well-known that there exists a long exact sequence \[\dots\rightarrow H_*(U)\rightarrow H_*(V)\rightarrow H_*(Cg)\rightarrow H_{*-1}(U)\rightarrow\dots
\rightarrow H_0(U)\rightarrow H_0(V)\rightarrow H_0(Cg)\rightarrow 0\]
(cf. \cite{Rotman}). The theorem is then a corollary of
Lemma~\ref{commutecone}. \epf

\section{Examples of calculation}
\label{example}

\begin{figure}
\begin{center}
\includegraphics[width=5cm]{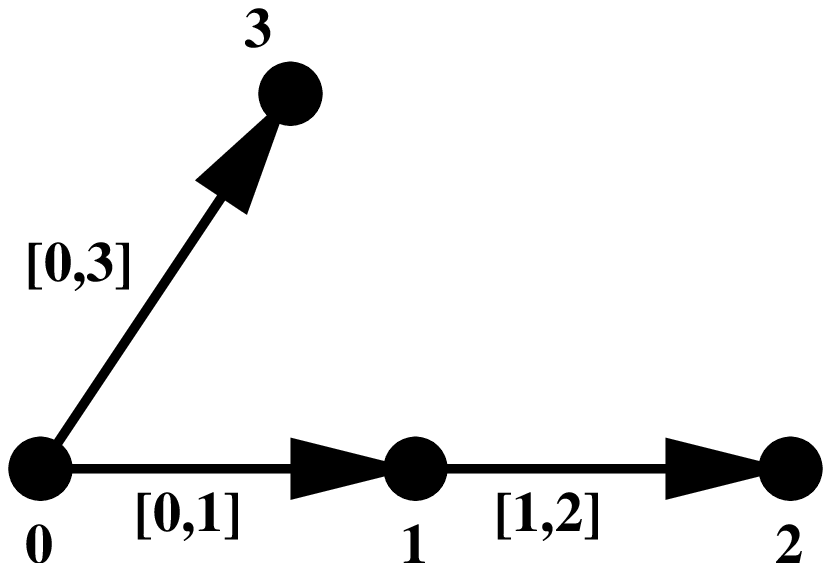}
\end{center}
\caption{$1$-dimensional branching}
\label{branch1}
\end{figure}

\begin{figure}
\begin{center}
\includegraphics[width=7cm]{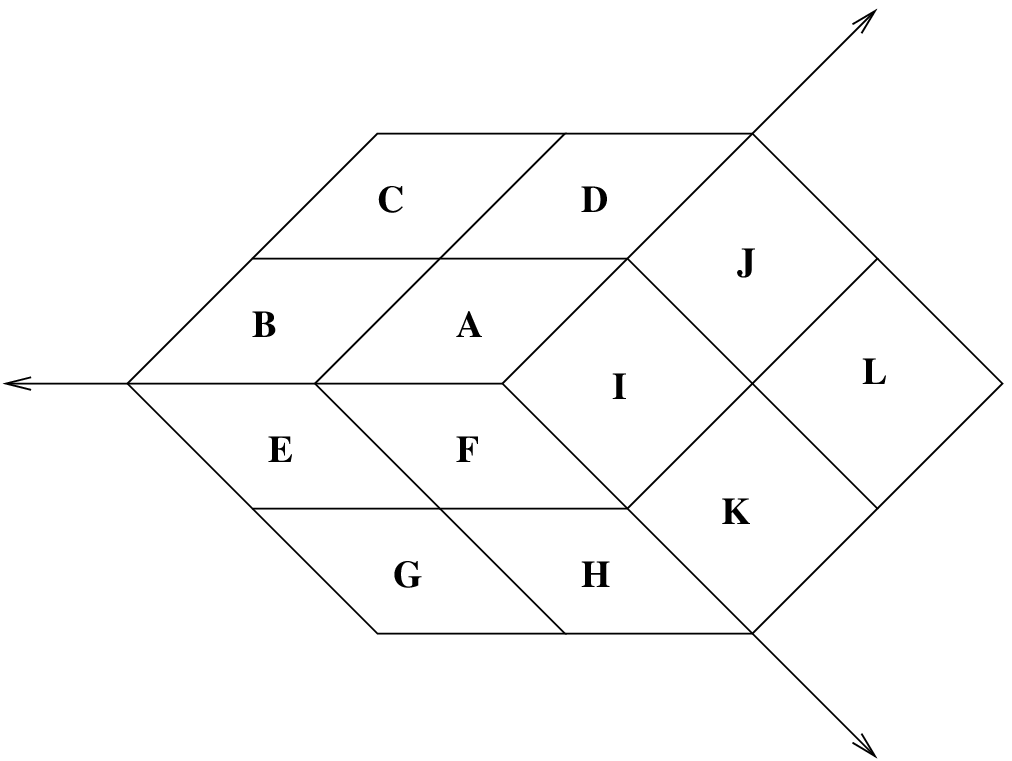}
\end{center}
\caption{$2$-dimensional branching}
\label{2branch}
\end{figure}

\begin{figure}
\begin{center}
\includegraphics[width=7cm]{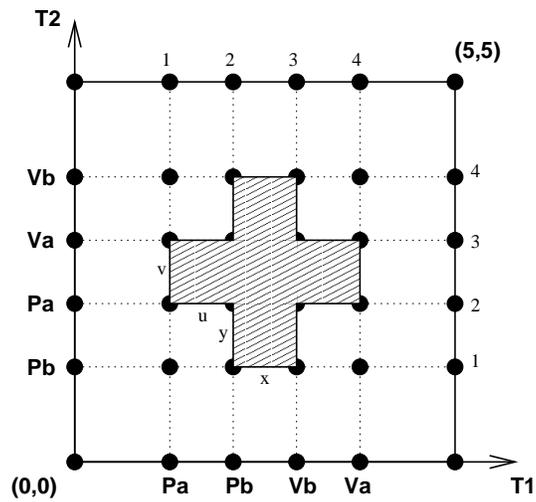}
\end{center}
\caption{The Swiss Flag Example}
\label{progress10}
\end{figure}

\bp\label{meme} 
If $X$ is a cofibrant flow, then the homotopy branching space
$\hop^-X$ and $\P^-X$ are homotopy equivalent. \ep

\bpf The functorial weak S-homotopy equivalence $Q(X)\longrightarrow X$ between 
cofibrant flows becomes a homotopy equivalence
$\P^-Q(X)\longrightarrow \P^-X$ of cofibrant topological spaces since
the functor $\P^-$ is a left Quillen functor. \epf

Since all examples given in this section are cofibrant flows, one can
then replace their homotopy branching space by their branching space.

\subsection{The directed segment}

By definition, the directed segment is the flow
\[\vI=\glob(\{[0,1]\}).\] One has $\P^-_0(\vI)=\{[0,1]\}$ and
$\P^-_1(\vI)=\varnothing$. And $H_n^-(\vI)=0$ for $n\geq 1$ and
$H_0^-(\vI)=\Z \{0,1\} / s(\P^-_0(\vI))$ is generated by the unique final
state of $\vI$.

\subsection{$1$-dimensional branching}

Consider the flow $X$ defined by $X^0=\{0,1,2,3\}$ and
$\P_{0,1}X=\{[0,1]\}$, $\P_{1,2}X=\{[1,2]\}$, $\P_{0,3}X=\{[0,3]\}$,
$\P_{0,2}X=\{[0,2]\}$ and $\P_{\alpha\beta}X=\varnothing$ otherwise
(cf. Figure~\ref{branch1}).

Then $\P^-_0X=\{[0,1],[0,3]\}$, $\P^-_1X=\{[1,2]\}$ and
$\P^-_2X=\P^-_3X=\varnothing$. One has $H_n^-(X)=0$ for $n\geq 2$, 
$H_1^-(X)=\Z$ (generated by $[0,3]-[0,1]$), and $H_0^-(X)=\Z\oplus\Z$ 
(generated by the final states $2$ and $3$).

\subsection{$2$-dimensional branching}

Let us consider now the case of Figure~\ref{2branch}. One has
$H_1^-=0$ and $H_n^-=0$ for $n\geq 2$. And $H_1^-=\Z$, the generating
branching being the one corresponding to the alternate sum
$(A)-(F)+(I)$. At last, $H_0^-=\Z\oplus \Z\oplus Z$, the generators
being the final states of the three squares $(C)$, $(G)$ and $(L)$. 
If $\alpha$ is the common initial state of $(A)$, $(F)$ and $(I)$, then 
$\P^-_\alpha=\mathbf{S}^1$.

\subsection{The Swiss Flag example}

Consider the discrete set \[SW^0=\{0,1,2,3,4,5\}\p \{0,1,2,3,4,5\}.\] 
Let
\beas
\mathcal{S}&=&\left\{((i,j),(i+1,j))\hbox{ for } (i,j)\in\{0,\dots,4\}\p \{0,\dots,5\}\right\}\\
&\cup& \left\{((i,j),(i,j+1))\hbox{ for } (i,j)\in\{0,\dots,5\}\p \{0,\dots,4\}\right\}\\
&\backslash & \left(
\{((2,2),(2,3)),((2,2),(3,2)), ((2,3),(3,3)),((3,2),(3,3))\}
\right)
\eeas
The flow $SW^1$ is obtained from $SW^0$ by attaching a copy of
$\glob(\mathbf{D}^0)$ to each pair $(x,y)\in \mathcal{S}$ with $x\in
SW^0$ identified with $0$ and $y\in SW^0$ identified with $1$.  The
flow $SW^2$ is obtained from $SW^1$ by attaching to each square
$((i,j),(i+1,j+1))$ except $(i,j)\in\{(2,1),(1,2),(2,2),(3,2),(2,3)\}$
a globular cell $\glob(\mathbf{D}^1)$ such that each execution path
$((i,j),(i+1,j),(i+1,j+1))$ and $((i,j),(i,j+1),(i+1,j+1))$ is
identified with one of the execution path of $\glob(\mathbf{S}^0)$
(there is not a unique choice to do that). Let $SW=SW^2$
(cf. Figure~\ref{progress10} where the bold dots represent the points
of the $0$-skeleton). The flow $SW$ represents the PV diagram of
Figure~\ref{progress10}.

The topological space $\P^-_\alpha$ is contractible for $\alpha\in
SW^0\backslash \{(1,2),(2,1),(5,5)\}$. And $\P^-_{(5,5)}=\varnothing$,
$\P^-_{(1,2)}=\{u,v\}$ and $\P^-_{(2,1)}=\{x,y\}$ with
$s(u)=s(v)=(1,2)$, $t(u)=(2,2)$, $t(v)=(1,3)$, $s(x)=s(y)=(2,1)$,
$t(x)=(3,1)$ and $t(y)=(2,2)$.

Then $H^-_0=\Z$ (generated by the final state $(5,5)$),
$H^-_1=\Z\oplus\Z$ (generated by $u-v$ and $x-y$). And $H_n^-=0$ for
any $n\geq 2$.

\section{Conclusion}

The branching homology is a dihomotopy invariant containing in
dimension $0$ the final states and in dimension $n\geq 1$ the
non-deterministic $n$-dimensional branching areas of non-constant
execution paths. The merging homology is a dihomotopy invariant
containing in dimension $0$ the initial states and in dimension $n\geq
1$ the non-deterministic $n$-dimensional merging areas of non-constant
execution paths. The non-deterministic branchings and mergings of
dimension $n\geq 2$ satisfies a long exact sequence which can be helpful for 
future applications or theoretical developments.

\appendix

\section{The case of mergings} 
\label{casemerging}

Some definitions and results about mergings are collected here,
almost without any comment or proof.

\bp \label{universalmp}
Let $X$ be a flow. There exists a topological space $\P^+X$ unique up
to homeomorphism and a continuous map $h^+:\P X\longrightarrow \P^+ X$
satisfying the following universal property:
\begin{enumerate}
\item For any $x$ and $y$ in $\P X$ such that $t(x)=s(y)$, the equality
$h^+(y)=h^+(x*y)$ holds.
\item Let $\phi:\P X\longrightarrow Y$ be a
continuous map such that for any $x$ and $y$ of $\P X$ such that
$t(x)=s(y)$, the equality $\phi(y)=\phi(x*y)$ holds. Then there exists a
unique continuous map $\overline{\phi}:\P^+X\longrightarrow Y$ such that
$\phi=\overline{\phi}\circ h^+$.
\end{enumerate}
Moreover, one has the homeomorphism
\[\P^+X\iso \bigsqcup_{\alpha\in X^0} \P^+_\alpha X\]
where $\P^+_\alpha X:=h^+\lp \bigsqcup_{\beta\in
X^0} \P^+_{\alpha,\beta} X\rp$. The mapping $X\mapsto \P^+X$
yields a functor $\P^+$ from $\dtop$ to $\top$. 
\ep

Loosely speaking, the merging space of a flow is the space of germs
of non-constant execution paths ending in the same way.

\bd \label{defplus}
Let $X$ be a flow. The topological space $\P^+X$ is called the {\rm
merging space} of the flow $X$. The functor $\P^+$ is called the {\rm
merging space functor}. \ed

Notice by that considering the opposite $X^{op}$ of a flow $X$ (by
interverting $s$ and $t$), then one obtains the following obvious
relation between $\P^-$ and $\P^+$ : $\P^+X=\P^-X^{op}$ and
$\P^-X=\P^+X^{op}$.

\bth  There exists a weak S-homotopy equivalence
of flows $f:X\longrightarrow Y$ such that the topological spaces
$\P^+X$ and $\P^+Y$ are not weakly homotopy equivalent. \eth

\bth The merging space functor $\P^+:\dtop\longrightarrow \top$ is a 
left Quillen functor. \eth

\bd
The {\rm homotopy merging space} $\hop^+ X$ of a flow $X$ is by
definition the topological space $\P^+Q(X)$.  If $\alpha\in X^0$, let 
$\hop^+_\alpha X=\P^+_\alpha X$. \ed

\begin{cor} 
Let $f:X\longrightarrow Y$ be a weak S-homotopy equivalence of flows.
Then $\hop^+ f:\hop^+ X\longrightarrow \hop^+ Y$ is a homotopy
equivalence between cofibrant topological spaces. \end{cor}

\bd Let $X$ be a flow. Then the $(n+1)$-th merging homology group
$H_{n+1}^+(X)$ is defined as the $n$-th homology group of the
augmented simplicial set $\mathcal{N}^+_*(X)$ defined as follows:
\begin{enumerate}
\item $\mathcal{N}^+_n(X)=\sing_n(\hop^+X)$ for $n\geq 0$
\item $\mathcal{N}^+_{-1}(X)=X^0$
\item the augmentation map $\epsilon:\sing_0(\hop^+X)\longrightarrow X^0$
is induced by the mapping $\gamma\mapsto s(\gamma)$ from $\hop^+X=\sing_0(\hop^+X)$
to $X^0$
\end{enumerate}
where $\sing(Z)$ denotes the singular simplicial nerve of a given topological space 
$Z$. In other terms, 
\begin{enumerate}
\item for $n\geq 1$, $H_{n+1}^+(X):=H_n(\hop^+X)$
\item  $H_1^+(X):=\ker(\epsilon)/\im\lp\partial:\mathcal{N}^+_1(X)\rightarrow
\mathcal{N}^+_0(X)\rp$
\item $H_0^+(X):=\Z(X^0)/\im(\epsilon)$.
\end{enumerate}
where $\partial$ is the simplicial differential map, where $\ker(f)$ is the kernel 
of $f$ and where $\im(f)$ is the kernel of $f$. 
\ed

\bp 
For any flow $X$, $H_0^+(X)$ is the free abelian group generated by the initial
states of $X$. 
\ep

\bp For any flow $X$, there exists a natural isomorphism of abelian groups 
\[H_{n+1}^+(X)\iso \bigoplus_{\alpha\in X^0} \widetilde{H}_n(\hop^+_\alpha X)\] 
for any $n\geq 0$.
\ep

\bp
Let $f:X\longrightarrow Y$ be a weak S-homotopy equivalence of flows. Then 
$\mathcal{N}^+(f):\mathcal{N}^+(X)\longrightarrow \mathcal{N}^+(Y)$ is a 
homotopy equivalence of augmented simplicial nerves. \ep

\begin{cor} 
Let $f:X\longrightarrow Y$ be a weak S-homotopy equivalence of
flows. Then $H_n^+(f):H_n^+(X)\longrightarrow H_n^+(Y)$ is an
isomorphism for any $n\geq 0$. \end{cor}

\bp Let $f:X\longrightarrow Y$ be a T-homotopy equivalence. Then for any
$n\geq 0$, the linear map $H_n^+(f):H_n^+(X)\longrightarrow H_n^+(Y)$
is an isomorphism.
\ep

\begin{cor} The merging homology is a dihomotopy invariant. \end{cor}

\begin{lem} 
The homotopy merging space of the terminal flow is
contractible. \end{lem}

\begin{lem} 
Let $f:X\longrightarrow Y$ be a morphism of flows. Let $Cf$ be the
cone of $f$.  Then the homotopy merging space $\hop^+ (Cf)$ of $Cf$
is homotopy equivalent to the cone $C (\hop^+f)$ of
$\hop^+f:\hop^+X\longrightarrow \hop^+Y$.
\end{lem}

\bth (Long exact sequence for higher dimensional mergings) 
For any morphism of flows $f:X\longrightarrow Y$, one has the long exact sequence
\beas
&& \dots \rightarrow H_{n}^+(X) \rightarrow H_{n}^+(Y) \rightarrow H_{n}^+(Cf)\rightarrow  \dots \\
&& \dots \rightarrow H_{3}^+(X) \rightarrow H_{3}^+(Y) \rightarrow H_{3}^+(Cf)\rightarrow \\
&& H_{2}^+(X) \rightarrow H_{2}^+(Y) \rightarrow H_{2}^+(Cf)\rightarrow \\
&& H_0(\hop^+X) \rightarrow H_0(\hop^+Y)\rightarrow  H_0(\hop^+ Cf)\rightarrow 0.
\eeas
\eth

We conclude this section by an additional remark about the Quillen
adjunctions induced by the functors $\P^-$ and $\P^+$.

\bth The Quillen adjunctions $\P^-:\dtop\rightleftarrows
\top:C^-$ and $\P^+:\dtop\rightleftarrows
\top:C^+$ together induce a Quillen adjunction
$\P^-\sqcup \P^+:\dtop\rightleftarrows
\top:C^-\p C^+$. \eth

\bpf Indeed, one has
\beas
\top(\P^-X\sqcup \P^+X,Z)&\iso & \top(\P^-X,Z)\p \top(\P^+X,Z)\\
&\iso & \dtop(X,C^-Z)\p \dtop(X,C^+Z)\\ &\iso & \dtop(X,C^-Z\p C^+Z)
\eeas
If $Z\longrightarrow T$ is a fibration of topological spaces, then
both $C^-Z\longrightarrow C^-T$ and $C^+Z\longrightarrow C^+T$ are
fibrations of flows by Theorem~\ref{rappel}. Since a product of fibrations
is a fibration, then  $C^-\p C^+$ is a right Quillen adjoint. And therefore
$\P^-\sqcup \P^+$ is a left Quillen adjoint.
\epf

None of the Quillen adjunctions $\P^-:\dtop\rightleftarrows
\top:C^-$, $\P^+:\dtop\rightleftarrows
\top:C^+$ and $\P^-\sqcup \P^+:\dtop\rightleftarrows
\top:C^-\p C^+$ gives rise to a Quillen equivalence. For obvious
reasons, the geometry of the branching space, the merging space or
both together cannot characterize a flow. Indeed, the information
about how branchings and mergings are related to one another is
missing.

\section{Branching space, merging space and S-homotopy}\label{correctforShomotopy}

The purpose of this section is to prove the:

\begin{prop} \label{correctforShomotopyprop}
Let $X$ and $Y$ be two S-homotopy equivalent flows
(cf. Definition~\ref{defS}) which are not necessarily cofibrant. Then
the topological spaces $\P^-X$ and $\P^-Y$ are homotopy
equivalent. \end{prop}

Proposition~\ref{correctforShomotopyprop} is already proved if $X$ and
$Y$ are both cofibrant: indeed since $\P^-:\dtop\longrightarrow\top$
is a left Quillen functor by Proposition~\ref{Pleft}, it preserves
weak equivalences between cofibrant objects.

Recall that two morphisms of flows $f,g:X\longrightarrow Y$ are
S-homotopy equivalent if and only if there exists a continuous map
$H:[0,1]\longrightarrow \tdtop(X,Y)$ such that $H(0)=f$ and $H(1)=g$
where the space $\tdtop(X,Y)$ is the set $\dtop(X,Y)$ equipped with
the Kelleyfication of the compact-open topology. In the same way, the
space $\ttop(U,V)$ denotes the set $\top(U,V)$ equipped with the
Kelleyfication of the compact-open topology. In particular, one has
the natural bijection of sets $\top(U\p V,W)\iso \top(U,\ttop(V,W))$
for any topological space $U$, $V$ and $W$.

We are going to need the category of non-contracting topological
$1$-categories.

\bd \cite{model3}
A {\rm non-contracting topological $1$-category} $X$ is a pair of
compactly generated topological spaces $(X^0,\P X)$ together with
continuous maps $s$, $t$ and $*$ satisfying the same properties as in
the definition of flow except that $X^0$ is not necessarily
discrete.  The corresponding category is denoted by $\cattopn$. \ed

\bp\label{tenseurn}\cite{model3}
Let $X$ and $Y$ be two objects of $\cattopn$. There exists a unique
structure of topological $1$-category
 $X\ot Y$ on the topological space  $X\p Y$ such that
\begin{enumerate}
\item $\left(X\ot Y\right)^0=X^0\p Y^0$ .
\item $\P\left(X\ot Y\right)= \left(\P X\p \P X\right)\sqcup \left(X^0 \p \P Y\right) \sqcup \left(\P X\p Y^0\right)$.
\item $s\left(x,y\right)=\left(s(x),s(y)\right)$, $t\left(x,y\right)=\left(t(x),t(y)\right)$, $\left(x,y\right)*\left(z,t\right)=\left(x*z,y*t\right)$.
\end{enumerate}
\ep

\bth\label{na}\cite{model3} The tensor product\index{tensor product} of
$\cattopn$ is a closed symmetric monoi\-dal structure, that is there
exists a bifunctor \[ \homcat:\cattopn\p \cattopn
\longrightarrow \cattopn\] contravariant with respect to the first
argument and covariant with respect to the second argument such that
one has the natural isomorphism of sets \[\cattopn\left(X\ot
Y,Z\right)\iso \cattopn\left(X,\homcat\left(Y,Z\right)\right)\] for
any topological $1$-categories $X$, $Y$ and $Z$. Moreover, one has the
natural homeomorphism
\[\lp \homcat\left(Y,Z\right)\rp^0\iso \tdtop(Y,Z).\]
\eth

With the tools above at hand, we can now prove the

\bth\label{continuite}
The functor $\P^-:\dtop\longrightarrow \top$ induces a natural
continuous map $(\P^-)_*:\tdtop(X,Y)\longrightarrow\ttop(\P^- X,\P^-
Y)$ for any flow $X$ and $Y$. \eth

\bpf The functor $\P^-:\dtop\longrightarrow \top$ yields
a set map \[\tdtop(X,Y)\longrightarrow\ttop(\P^- X,\P^- Y).\] One has
to prove that this set map is continuous.

By Yoneda's lemma, one has an isomorphism between the set
\[\nat\lp\top\lp-, \tdtop(X,Y)\rp,\top\lp-, \ttop(\P^-X,\P^-Y)\rp\rp\]
and the set \[\top\lp \tdtop(X,Y),\ttop(\P^- X,\P^- Y)\rp\] where
$\nat(F,G)$ denotes the set of natural transformations from a functor
$F$ to  another functor $G$.

Let $U$ be a topological space. Then $U$ can be viewed as a
non-contracting topological $1$-category if $U$ is identified with its
$0$-skeleton. Then
\beas
\top\lp U, \tdtop(X,Y)\rp &\iso & \top\lp U, \lp \homcat\left(X,Y\right)\rp^0\rp\\
&\iso & \cattop\lp U,\homcat(X,Y)\rp \\
& \iso & \cattop\lp U\otimes X,Y\rp.
\eeas
Let $\P^-:\cattop\longrightarrow \top$ be the functor defined as follows:
if $X$ is an object of $\cattop$, then the topological space $\P^-X$ is the
quotient of the topological space $\P X$ by the topological closure of
the smallest equivalence relation identifying $x$ and $x*y$ for any
$x,y\in\P X$ such that $t(x)=s(y)$. Clearly, one has the commutative
diagram of functors
\[
\xymatrix{
\dtop \fr{\P^-}\fd{} & \top \fd{=}\\
\cattop \fr{\P^-} & \top
}
\]
where the functor $\dtop\longrightarrow \cattop$ is the canonical embedding.

The non-contract\-ing topological $1$-category $U\otimes X$ looks as
follows: the $0$-skeleton is equal to $U\p X^0$ and the path
space is equal to $U\p \P X$ with the composition law characterized by
$s(u,x)=(u,sx)$, $t(u,x)=(u,tx)$ and $(u,x)*(u,y)=(u,x*y)$. Therefore
there exists a natural homeomorphism $\P^-(U\otimes X)\iso U\p
\P^-X$. So the functor $\P^-:\cattop\longrightarrow \top$ induces a
set map
\[\cattop\lp U\otimes X,Y\rp\longrightarrow \top\lp U\p \P^-X,\P^-Y\rp\]
Since $\top\lp U\p \P^-X,\P^-Y\rp\iso \top\lp U,\ttop\lp
\P^-X,\P^-Y\rp\rp$, one obtains by composition a natural set map
\[\top\lp U, \tdtop(X,Y)\rp\longrightarrow
\top\lp
U,\ttop\lp \P^-X,\P^-Y\rp\rp\] which by Yoneda's lemma provides a
continuous map \[\tdtop(X,Y)\longrightarrow \ttop\lp \P^-X,\P^-Y\rp\]
whose underlying set map is exactly the set map
$\dtop(X,Y)\longrightarrow\top\lp \P^-X,\P^-Y\rp$ induced by the
functor $\P^-:\dtop\longrightarrow \top$.
\epf

\begin{cor} Let $f$ and $g$ be two S-homotopy equivalent morphisms of flows
from $X$ to $Y$. Then the continuous maps $\P^-f$ and $\P^-g$ from
$\P^-X$ to $\P^-Y$ are homotopic. \end{cor}

\bpf Let $H$ be an element of $\top([0,1],\tdtop(X,Y))$ such that $H(0)=f$ and
$H(1)=g$. Then $(\P^-)_*(H)\in \top([0,1],\ttop(\P^-X,\P^-Y))$ yields
an homotopy from $\P^-f$ to $\P^-g$.
\epf

\begin{cor}
Let $X$ and $Y$ be two S-homotopy equivalent flows. Then the
topological spaces $\P^-X$ and $\P^-Y$ are homotopy equivalent.
\end{cor}

Of course, the same theorem holds for the merging space functor:

\begin{cor}
Let $X$ and $Y$ be two S-homotopy equivalent flows. Then the
topological spaces $\P^+X$ and $\P^+Y$ are homotopy equivalent.
\end{cor}

\end{document}